\documentclass[10pt]{amsart}
\usepackage{latexsym,amsmath,amssymb}

 \renewcommand{\epsilon}{\varepsilon}
 \newcommand{\newsection}[1]
  {\section{#1}\setcounter{theorem}{0} \setcounter{equation}{0}\par\noindent}

   \newtheorem{theorem}{Theorem}[section]
   \newtheorem{lemma}[theorem]{Lemma}
 \newtheorem{corr}[theorem]{Corollary}
 
 \newtheorem{proposition}[theorem]{Proposition}
 \newtheorem{deff}[theorem]{Definition}
 \newtheorem{remark}[theorem]{Remark}

  \numberwithin{equation}{section}

 \newcommand{\bth}{\begin{theorem}}
 \newcommand{\ble}{\begin{lemma}}
 \newcommand{\bcor}{\begin{corr}}
 \newcommand{\bdeff}{\begin{deff}}
 \newcommand{\bprop}{\begin{proposition}}
 \def\be{\begin{equation}}
\def\ee{\end{equation}}
\def\bt{\begin{theorem}}
\def\et{\end{theorem}}
\def\ba{\begin{array}}
\def\ea{\end{array}}
\def\bl{\begin{lemma}}
\def\el{\end{lemma}}

 \newcommand{\ele}{\end{lemma}}
 \newcommand{\ecor}{\end{corr}}
 \newcommand{\edeff}{\end{deff}}
 
 \newcommand{\eprop}{\end{proposition}}

 \renewcommand{\Pi}{\varPi}

 \renewcommand{\epsilon}{\varepsilon}

\title[Singularities of compressible Euler equations]
{  Singularities of solutions to compressible Euler equations with
vacuum}

\date{\today}

\begin{document}
\maketitle

\centerline{\author{Zhen
 Lei\footnote{School of Mathematical Sciences; LMNS and Shanghai
 Key Laboratory for Contemporary Applied Mathematics, Fudan University, Shanghai 200433, P. R. China. {\it Email:
 leizhn@gmail.com}}}\ and Yi Du\footnote{School of Mathematical Sciences, South China Normal University, Guangzhou 510631, P. R. China.
 {\it Email: duyidy@gmail.com}}
 \and Qingtian Zhang\footnote{Department of Mathematics, Penn State University, State College, PA 16801, USA. {\it Email:
 zhang\_q@math.psu.edu}}}

\begin{abstract}
Presented are two results on the formation of finite time
singularities of solutions to the compressible Euler equations in
two and three space dimensions for isentropic, polytropic, ideal
fluid flows. The initial velocity is assumed to be symmetric and
the initial sound speed is required to vanish at the origin. They
are smooth in Sobolev space $H^3$, but not required to have a
compact support. It is shown that the $H^3$ norm of the velocity
field and the sound speed will blow up in a finite time.
\end{abstract}

{\bf Keywords: } Finite time singularities, compressible Euler
equations, vacuum.

\newsection{Introduction}

Euler equation is one of the most fundamental equations in fluid
dynamics. Many interesting fluid dynamic phenomena can be
described by the Euler equation (see, for instance,
\cite{Lions98}). Recently, singularity formation in fluid
mechanics has attracted the attention of a number of researchers,
see, for instance,
\cite{CCFGG11-b,Constantin86,CFM96,CS12,HL12,HL09a} and two recent
review articles \cite{BT07,Constantin07}. For compressible Euler
equations of the motion of polytropic ideal fluid flows, Sideris
(see \cite{Sideris}) proved, under various settings, several very
interesting results on the formation of finite time singularities
to solutions whose initial velocity field has a compact support
and initial density is strictly positive and is equal to a
positive constant outside the support of the initial velocity
field. It is very interesting to investigate the long time
behavior of solutions to compressible Euler equations with initial
data containing vacuum states, as has been pointed out in
\cite{Mas12}.

In this short article, we prove that solutions to compressible
Euler equations will develop finite time singularities for
radially symmetric initial data whose initial velocity field has
no compact support and initial density contains vacuum states. In
particular, it is shown that the $H^3$ norm of the velocity field
and the sound speed will blow up in a finite time.

To state our theorems, let us  begin with the compressible Euler
equations for isentropic, polytropic, ideal fluid flows:
\begin{equation}\label{Euler}
\begin{cases}
\rho_t+\nabla\cdot(\rho u)=0,
\\
(\rho u)_t+\nabla\cdot (\rho u\otimes u)+\nabla p=0,
\end{cases}
\end{equation}
where $\rho(t, \cdot): \mathbb{R}^n \rightarrow \mathbb{R}$ is the
scalar mass density, $u(t, \cdot): \mathbb{R}^n \rightarrow
\mathbb{R}^n$ is the velocity field,  $p(t, x)$ is the pressure
which is given by the equation of state
\begin{equation}\label{ES}
p = A\rho^\gamma.
\end{equation}
Here $A > 0$ is an entropy constant, $\gamma$ is the adiabatic
index. For polytropic gases, one has $1 < \gamma \leq
\frac{5}{3}$.

The compressible Euler equations \eqref{Euler} are imposed on the
following initial data
\begin{equation}\label{data}
\rho(0, x) = \rho_0(r),\quad u(0, x) = \frac{x}{r}v_0(r).
\end{equation}
Here and in what follows, we will use
\begin{equation}\nonumber
r = |x|,\quad x \in \mathbb{R}^n,
\end{equation}
for notational convenience. Denote the sound speed $c(\rho)$ by
\begin{equation}\label{SS}
c(\rho) = \sqrt{\frac{\partial p(\rho)}{\partial \rho}} =
\sqrt{A\gamma}\rho^{\frac{\gamma - 1}{2}}.
\end{equation}

Our first result is on the formation of singularities of solutions
to the compressible Euler equations in three dimensions:
\begin{theorem}\label{Threed}
Assume that $\gamma > 1$, $\rho_0 \geq 0$ and $(c_0, u_0) \in
H^3(\mathbb{R}^3)$ with $c_0 = c(\rho_0)$. Moreover, assume that
$\rho_0$ and $v_0$ satisfy
\begin{equation}\label{cond-1}
\rho_0(0) = 0,
\end{equation}
\begin{equation}\label{cond-2}
\int_{\mathbb{R}^3} \rho_0(r)dx > 0,
\end{equation}
and
\begin{equation}\label{cond-3}
- \int_{\mathbb{R}^3}\frac{(1 + r)\rho_0v_0}{r^2e^r}dx \geq
\sqrt{\frac{A}{(\gamma + 1)
  (4\pi)^{\gamma-1}}}\big(
  \int_{\mathbb{R}^3}\frac{\rho_0}{re^{r}}
  dx\big)^{\frac{\gamma + 1}{2}}.
\end{equation}
Then the solution $(\rho, u)$ to the compressible Euler equations
\eqref{Euler} with the initial data \eqref{data} will develop
finite time singularities.
\end{theorem}

Our second result is on the formation of singularities of
solutions to the compressible Euler equations in two dimensions:
\begin{theorem}\label{Twod}
Assume that $\gamma > 1$, $\rho_0 \geq 0$ and $(c_0, u_0) \in
H^3(\mathbb{R}^2)$ with $c_0 = c(\rho_0)$. Let $K_0(r)$ be the
modified Bessel function
\begin{equation}\label{BF}
K_0(r)=\int_0^\infty e^{-r\cosh t}dt.
\end{equation}
If the initial data $\rho_0$ and $v_0$ satisfy \eqref{cond-1},
\begin{equation}\label{cond-5}
\int_{\mathbb{R}^2}\rho_0(r)dx > 0,
\end{equation}
and
\begin{equation}\label{cond-6}
\int_{\mathbb{R}^2}\rho_0(r)v_0(r)K_0^\prime(r)dx \geq
\sqrt{\frac{A}{\gamma + 1}}\frac{\big(
  \int_{\mathbb{R}^2}\rho_0(r)K_0(r)
  dx\big)^{\frac{\gamma + 1}{2}}}{
  \big(\int_{\mathbb{R}^2}K_0(r)dx\big)^{\frac{\gamma-1}{2}}},
\end{equation}
then the solution $(\rho, u)$ to the compressible Euler equations
\eqref{Euler} with the initial data \eqref{data} will develop
finite time singularities.
\end{theorem}

\begin{remark}\label{rem}
For a polytropic ideal gas, the adiabatic index $\gamma \in
(1,\frac{5}{3}]$ and hence $\frac{2}{\gamma - 1} \geq 3$. Due to
the expression of the sound speed in \eqref{SS}, it is easy to get
\begin{equation}\nonumber
\rho_0 = \big(A\gamma \big)^{- \frac{1}{\gamma -
1}}c_0^{\frac{2}{\gamma - 1}}.
\end{equation}
Consequently, one also has $\rho_0 \in H^3(\mathbb{R}^n)$ under
the condition that $c_0 \in H^3(\mathbb{R}^n)$, which in turn
implies that $\rho(t, \cdot) \in H^2(\mathbb{R}^n)$ as long as the
solution is smooth. For $\gamma > \frac{5}{3}$, local
well-posedness theory in Theorem \ref{thm-local} implies that
$\rho \in C([0, T) \times \mathbb{R}^n)$.
\end{remark}

\begin{remark}
For the non-isentropic, polytropic, ideal gases, the entropy $S$
is transported by the flows. It is easy to verify that our proofs
of the blow up parts in Theorem \ref{Threed} and Theorem
\ref{Twod} are still true. However, the local well-posedness of
the coupled system with initial vacuum is much more complicated
(see, for instance, \cite{Mas12}). We do not pursue this issue in
this short article.
\end{remark}

There are several ingredients in the proofs of the above theorems.
The first one is to write the compressible Euler equations
\eqref{Euler} as a quasi-linear wave type equation in terms of
$\rho$ with inhomogeneous terms involving $u$:
\begin{equation}\nonumber
\rho_{tt}-\Delta p =\nabla\cdot[\nabla\cdot(\rho u\otimes u)].
\end{equation}
This is in fact what Sideris did in \cite{Sideris}. However, we
will treat this equation in a very different manner from that in
\cite{Sideris} due to the facts that $\rho(r) \rightarrow 0$ as $r
\rightarrow \infty$ and the initial velocity field $u_0$ has no
compact support. We will choose the modified Bessel function
$K_0(r)$ in 2D case and $\frac{1}{re^{r}}$ in 3D case as test
functions for the above quasi-linear wave type equation,
respectively, to explore the nonlinear structure of the pressure
$p$ as a function of $\rho$, which eventually corresponds to the
formation of finite time singularities. The second one is making
use of the symmetric structure of the solutions, which results in
a good sign for the inhomogeneous term
$\nabla\cdot[\nabla\cdot(\rho u\otimes u)]$ when taking the inner
product of the above wave type equation with test functions.
Moreover, the symmetric structure of the solutions will also be
used to eliminate boundary terms on the artificial boundary $r =
0$.

The remaining part of this paper is simply organized as follows:
In section 2 we will prove Theorem \ref{Threed}. Then in section 3
we present the proof of Theorem \ref{Twod}.


\section{Singularities of compressible Euler equations in 3D}

In this section we will prove Theorem \ref{Threed}. Before that,
let us recall the local well-posedness of compressible Euler
equations (for example, see \cite{Kato75}). To be precise, let us
formulate it as a theorem for our use.

\begin{theorem}\label{thm-local}
Assume that $\gamma > 1$, $\rho_0 \geq 0$ and $(c_0, u_0) \in
H^3(\mathbb{R}^n)$ with $c_0 = c(\rho_0)$ and $n = 2,\ 3$. Then
there exists a unique solution $(\rho, u)$ to the compressible
Euler equations \eqref{Euler} with initial data \eqref{data}  on
some time interval $[0, T)$, which satisfies
\begin{equation}\nonumber
c(\rho),\ u \in C([0, T), H^3(\mathbb{R}^n)) \cap C^1([0, T),
H^2(\mathbb{R}^n)) \cap C^2([0, T), H^1(\mathbb{R}^n)),
\end{equation}
and
\begin{equation}\label{CoD}
\rho \in C([0, T) \times \mathbb{R}^n).
\end{equation}
If $1 < \gamma \leq \frac{5}{3}$, then
\begin{equation}\label{3-1}
\rho \in C([0, T), H^2(\mathbb{R}^n)) \cap C^1([0, T),
H^1(\mathbb{R}^n)).
\end{equation}
If $\rho_0$ and $u_0$ are radially symmetric and \eqref{cond-1} is
satisfied, then
\begin{equation}\label{3-2}
\rho(t, 0) \equiv 0,\ u(t, 0)  \equiv 0,\quad u(t, x) =
\frac{x}{r}v(t, r).
\end{equation}
\end{theorem}

The proof of the local well-posedness in the above theorem is
based on the fact that the compressible Euler equations can be
written as a symmetric hyperbolic system in terms of $(c, u)$:
\begin{equation}\nonumber
\begin{cases}
c_t + u\cdot\nabla c + \frac{2c}{\gamma - 1}\nabla\cdot u = 0,\\
u_t + u\cdot\nabla u + \frac{\gamma c}{\gamma - 1}\nabla c = 0.
\end{cases}
\end{equation}
See \cite{Kato75} for more details. The fact \eqref{3-1} is due to
Remark \ref{rem} and the equation of mass conservation in the
 original compressible Euler equations \eqref{Euler}. The fact
$u(t, 0) \equiv 0$ is in fact a universal identity for smooth
symmetric vector. To see $\rho(t, 0) \equiv 0$, we use $u(t, 0)
\equiv 0$ and the equation of mass conservation to get
\begin{equation}\nonumber
\rho(t, 0) = \rho_0(0)e^{- \int_0^t\nabla\cdot u(s, 0)ds}.
\end{equation}

 It is ready to present the proof of Theorem \ref{Threed}.

\begin{proof}
We prove Theorem \ref{Threed} by contradiction. Suppose that the
solution $(c, u) \in H^3(\mathbb{R}^3)$ for all time $t \geq 0$
and $T = \infty$ in Theorem \ref{thm-local}. We will derive that
the density blows up in the ball centered at the origin with an
arbitrary small radius $r_0 > 0$ in a finite time, which
contradicts with \eqref{CoD}.

Applying the time derivative to the first equation of the
compressible Euler system $\eqref{Euler}$,
 we have (in the sense of distribution)
 \begin{equation}\label{3-3}
 \rho_{tt}=-\nabla\cdot(\rho u)_t=\Delta p +\nabla\cdot[\nabla\cdot(\rho u\otimes u)].
 \end{equation}
Taking the $L^2$ inner product of the above equation with the test
function $\frac{1}{re^{r}}$, one has
\begin{align}\nonumber
\frac{d^2}{dt^2}\int_{\mathbb{R}^3}\frac{\rho}{re^{r}}
dx&=\int_{\mathbb{R}^3}\Delta p\frac{1}{re^{r}} dx
+\int_{\mathbb{R}^3}\nabla\cdot[\nabla\cdot(\rho u\otimes u)]
\frac{1}{re^{r}} dx.
\end{align}
Let us first compute that
\begin{eqnarray}\nonumber
&&\int_{\mathbb{R}^3}\Delta p\frac{1}{re^{r}}
  dx\\\nonumber
&&=\int_{\mathbb{R}^3}p\Delta \frac{1}{re^{r}}
  dx - \lim_{\epsilon\rightarrow0}\int_{|x|=\epsilon}\big(\frac{\partial p}
  {\partial r}\frac{e^{-r}}{r} -\frac{d(\frac{1}{re^{r}})}{d r}p\big)ds\\\nonumber
&&=\int_{\mathbb{R}^3}p
  \frac{1}{re^{r}}dx+\lim_{\epsilon\rightarrow0}
  \int_{|x|=\epsilon}(\frac{1}{re^{r}})^\prime p ds.
\end{eqnarray}
Using the equation of state \eqref{ES} and the fact that
$\rho(t,0)\equiv0$, we have
\begin{eqnarray}\nonumber
&&A\int_{\mathbb{R}^3}\Delta \rho^\gamma\frac{1}{re^{r}} dx =
  A\int_{\mathbb{R}^3}\rho^\gamma \frac{1}{re^{r}}dx\\\nonumber
&&\geq \frac{A}{(4\pi)^{\gamma-1}}\big(
  \int_{\mathbb{R}^3}\frac{\rho}{re^{r}}
  dx\big)^\gamma.
\end{eqnarray}
On the other hand,  noting \eqref{3-2}, one has
\begin{align}\nonumber
&\int_{\mathbb{R}^3}\frac{1}{re^{r}}\nabla\cdot[\nabla\cdot(\rho
  u\otimes u)] dx\\\nonumber
=& \int_{\mathbb{R}^3}
\partial_i\partial_j(\frac{1}{re^{r}})(\rho u_iu_j)dx\\\nonumber
=&\int_{\mathbb{R}^3}\rho v^2
\frac{x_j}{r}\partial_r[\frac{x_j}{r}(\frac{1}{re^{r}})^\prime]
dx\\\nonumber =&\int_{\mathbb{R}^3}\rho v^2
(\frac{1}{re^{r}})^{\prime\prime}dx\\\nonumber
=&\int_{\mathbb{R}^3}\rho
v^2[\frac1r+\frac2{r^2}+\frac2{r^3}]e^{-r}dx>0.
\end{align}
Here and in what follows, we use Einstein's convention for
summation over repeated indices. Consequently,
\begin{equation}\label{3-4}
\frac{d^2}{dt^2}F(t) \geq \frac{A}{(4\pi)^{\gamma-1}}\big(
F(t)\big)^\gamma,
\end{equation}
where
\begin{equation}\nonumber
F(t) = \int_{\mathbb{R}^3}\frac{\rho(t, r)}{re^{r}} dx.
\end{equation}

Using the equation of conservation of mass and integration by
parts, and noting \eqref{cond-3}, one has
\begin{equation}\label{3-5}
F^\prime(0) = \frac{d}{dt}\int_{\mathbb{R}^3}\frac{\rho}{re^{r}}
dx\Big|_{t = 0} = - \int_{\mathbb{R}^3}\frac{(1 + r)\rho_0
v_0}{r^2e^{r}}dx>0.
\end{equation}
The combination of \eqref{3-4} and \eqref{3-5} gives that
\begin{eqnarray}\nonumber
&&F^\prime(t) = F^\prime(0) +  \int_0^t \frac{d^2}{ds^2}F(s) ds
\geq F^\prime(0) > 0.
\end{eqnarray}
Consequently, one can multiply the both sides of \eqref{3-4} by
$\frac{d}{dt}F(t)$ to get
\begin{eqnarray}\nonumber
&&\big(F^\prime(t)\big)^2\geq \frac{A}{(\gamma + 1)
  (4\pi)^{\gamma-1}}\big(F(t)\big)^{\gamma + 1}\\\nonumber
&&\quad + \Big(\frac{d}{dt}\int_{\mathbb{R}^3}\frac{\rho}{re^{r}}
  dx\Big)^2\Big|_{t = 0} - \frac{A}{(\gamma + 1)
  (4\pi)^{\gamma-1}}\big(
  \int_{\mathbb{R}^3}\frac{\rho_0}{re^{r}}
  dx\big)^{\gamma + 1}.
\end{eqnarray}
Noting \eqref{cond-3} and \eqref{3-5}, and by denoting
\begin{eqnarray}\nonumber
C_0 = \sqrt{\frac{A}{(\gamma + 1) (4\pi)^{\gamma-1}}},
\end{eqnarray}
one has
\begin{eqnarray}\nonumber
F^\prime(t) \geq C_0\big(F(t)\big)^{\frac{\gamma + 1}{2}},
\end{eqnarray}
which implies that
\begin{equation}\nonumber
\int_{\mathbb{R}^3}\frac{\rho}{re^{r}} dx \geq \big(F(0)^{-
\frac{\gamma - 1}{2}} - \frac{\gamma - 1}{2}C_0t\big)^{-
\frac{2}{\gamma - 1}}.
\end{equation}
Noting \eqref{cond-2}, one has
\begin{equation}\nonumber
\int_{\mathbb{R}^3} \frac{\rho_0}{re^r}dx > 0,
\end{equation}
By using the mass conservation, one has
\begin{equation}\nonumber
F(t) = \int_{\mathbb{R}^3}\frac{\rho}{re^{r}} dx \leq
\int_{B_{r_0}}\frac{\rho}{re^{r}} dx +
\frac{1}{r_0}\int_{\mathbb{R}^3}\rho_0 dx,
\end{equation}
for any given $r_0 > 0$, here $B_{r_0}$ is the three-dimensional
ball centered at the origin with radius $r_0$. Consequently, one
concludes that $\int_{r \leq r_0}\rho(t, r)rdr$ can not be bounded
as $t \rightarrow \frac{2F(0)^{- \frac{\gamma - 1}{2}}}{(\gamma -
1)C_0}$, which in turn implies that $\rho(t, r)$ will blow up for
$r \leq r_0$ as $t \rightarrow \frac{2F(0)^{-\frac{\gamma -
1}{2}}}{(\gamma - 1)C_0}$. By the assumption $c \in L^\infty([0,
\infty), H^3(\mathbb{R}^3)$ and Sobolev imbedding theorem, one has
$c \in L^\infty([0, \infty)\times \mathbb{R}^3)$ and hence $\rho
\in L^\infty([0, \infty)\times \mathbb{R}^3)$. We arrive at a
contradiction. So the $H^3$ norm of $(c, u)$ can not be bounded
before the time $\frac{2F(0)^{-\frac{\gamma - 1}{2}}}{(\gamma -
1)C_0}$. The proof of Theorem \ref{Threed} is completed.
\end{proof}

\section{Singularities of compressible Euler equations in 2D}

In this section we prove Theorem \ref{Twod}. First of all, let us
recall that the modified Bessel function $K_0(r)$ and its
derivative $K_0^\prime(r)$ decay sufficiently fast as $r
\rightarrow \infty$. In fact, for any given $k > 1$, one has
\begin{eqnarray}\nonumber
&&\sup_{0 < r < \infty}r^kK_0(r) = \sup_{0 < r <
  \infty}r^k\int_0^\infty e^{-r\cosh t}dt\\\nonumber
&&= \sup_{0 < r < \infty} r^k\int_0^1 e^{-\frac{re^t}{2}}
  e^{-\frac{re^{-t}}{2}}dt + \sup_{0 < r < \infty}r^k
  \int_1^\infty e^{-\frac{re^t}{2}}e^{-\frac{re^{-t}}{2}}dt\\\nonumber
&&\leq \sup_{0 < r < \infty}r^ke^{- \frac{r}{2e}}
  \int_0^1e^{-\frac{re^t}{2}}dt + \sup_{0 < r < \infty}
  \int_1^\infty r^ke^{-\frac{re^t}{2}}(\frac{r}{2}e^t)^k(\frac{2}{r}e^{-t})^kdt\\\nonumber
&&\leq \sup_{0 < r < \infty}r^ke^{- \frac{r}{2e}}
  + \sup_{0 < r < \infty}
  2^k\int_1^\infty e^{-kt}dt\sup_{0 < s < \infty}e^{-s}s^k\\\nonumber
&&< \infty
\end{eqnarray}
and
\begin{eqnarray}\nonumber
&&\sup_{0 < r < \infty}r^k|K_0^\prime(r)| = \sup_{0 < r <
  \infty}r^k\int_0^\infty e^{-r\cosh t}\cosh tdt\\\nonumber
&&= \sup_{0 < r < \infty} 2r^k\int_0^1 e^{-\frac{re^t}{2}}
  e^{-\frac{re^{-t}}{2}}dt + \sup_{0 < r < \infty}r^k
  \int_1^\infty e^{-\frac{re^t}{2}}e^{-\frac{re^{-t}}{2}}e^tdt\\\nonumber
&&\leq \sup_{0 < r < \infty}2r^ke^{- \frac{r}{2e}}
  + \sup_{0 < r < \infty}
  2^k\int_1^\infty e^{-(k-1)t}dt\sup_{0 < s < \infty}e^{-s}s^k\\\nonumber
&&< \infty.
\end{eqnarray}
On the other hand, for $r > 0$, we also have
\begin{eqnarray}\nonumber
&&K_0(r)\leq \int_0^1 e^{-r\cosh t}dt + \int_1^\infty
  e^{-\frac{re^t}{2}}e^{-\frac{re^{-t}}{2}}dt\\\nonumber
&&\leq 1 + \int_1^\infty e^{-\frac{re^t}{2}}dt \leq 1 +
\frac{2}{r},
\end{eqnarray}
and
\begin{eqnarray}\nonumber
&&|K_0^\prime(r)| \leq \int_0^1 e^{-r\cosh t}e^tdt +
  \int_1^\infty e^{-\frac{re^t}{2}}e^{-\frac{re^{-
  t}}{2}}e^tdt\\\nonumber
&&\leq 2 + \int_1^\infty e^{-\frac{re^t}{2}}e^tdt\\\nonumber
&&\leq 2 + \int_1^\infty 2\big(\frac{re^t}{2}\big)^{-2}e^tdt \leq
2 + \frac{1}{2r^2}.
\end{eqnarray}

We in fact have proved the following lemma:
\begin{lemma}\label{lem}
The modified Bessel function $K_0(r) = \int_0^\infty e^{-r\cosh
t}dt$ satisfies
\begin{equation}\nonumber
\begin{cases}
K_0(r) \leq \frac{3}{r}, |K_0^\prime(r)| \leq \frac{1}{r^2},\quad
0 < r < \frac{1}{2},\\
K_0(r) \leq \frac{C_k}{r^k}, |K_0^\prime(r)| \leq
\frac{C_k}{r^k},\quad r > 1,
\end{cases}
\end{equation}
for constants $C_k$ depending only on $k > 1$.
\end{lemma}

Now it is ready to present the proof of Theorem \ref{Twod}.

\begin{proof}
We prove Theorem \ref{Twod} by contradiction. Suppose that the
solution $(c, u) \in H^3(\mathbb{R}^3)$ for all time $t \geq 0$
and $T = \infty$ in Theorem \ref{thm-local}. We will derive that
the density blows up in the ball centered at the origin with an
arbitrary small radius $r_0 > 0$ in a finite time, which
contradicts with \eqref{CoD}.

Taking the $L^2$ inner product of the above equation with the test
function $K_0(r)$, one has
\begin{align}\nonumber
\frac{d^2}{dt^2}\int_{\mathbb{R}^2}\rho K_0(r)
dx&=\int_{\mathbb{R}^2}\Delta pK_0(r)dx
+\int_{\mathbb{R}^2}\nabla\cdot[\nabla\cdot(\rho u\otimes u)]
K_0(r)dx.
\end{align}
Using a similar argument as in section 2 and noting the decay
properties of the Bessel function $K_0(r)$ in Lemma \ref{lem}, we
have
\begin{eqnarray}\nonumber
\int_{\mathbb{R}^2}\Delta pK_0(r)
  dx = \int_{\mathbb{R}^2}p
  \Delta K_0(r)dx.
\end{eqnarray}
Noting that the modified Bessel function $K_0(r)$ satisfies
\begin{equation}\label{4-1}
K^{\prime\prime}_0 + \frac{1}{r}K_0^\prime = K_0,
\end{equation}
one has
\begin{eqnarray}\nonumber
\int_{\mathbb{R}^2}\Delta p K_0(r)
  dx = A\int_{\mathbb{R}^2}\rho^\gamma K_0(r)dx\\\nonumber
\geq
\frac{A}{\big(\int_{\mathbb{R}^2}K_0(r)dx\big)^{\gamma-1}}\big(
  \int_{\mathbb{R}^2}\rho K_0(r)
  dx\big)^\gamma.
\end{eqnarray}
On the other hand,  noting \eqref{3-2} and using \eqref{4-1} and
Lemma \ref{lem}, one has
\begin{align}\nonumber
&\int_{\mathbb{R}^2}K_0(r)\nabla\cdot[\nabla\cdot(\rho
  u\otimes u)] dx\\\nonumber
=&\int_{\mathbb{R}^2}\rho v^2 K_0^{\prime\prime}(r)dx\\\nonumber
=&\int_{\mathbb{R}^2}\rho v^2\big(K_0(r) -
\frac{1}{r}K_0^\prime(r)\big)dx>0.
\end{align}
Here we also used the fact that $K_0^\prime(r) < 0$ by the
expression of the modified Bessel function $K_0(r)$. Consequently,
we have
\begin{equation}\label{4-4}
\frac{d^2}{dt^2}G(t) \geq
\frac{AG(t)^\gamma}{\big(\int_{\mathbb{R}^2}K_0(r)dx\big)^{\gamma-1}}.
\end{equation}
where
\begin{equation}\nonumber
G(t) = \int_{\mathbb{R}^2}\rho(t, r) K_0(r)dx.
\end{equation}

Using the equation of conservation of mass and integration by
parts, and noting \eqref{cond-5}, one has
\begin{equation}\label{4-5}
G^\prime(0) = \frac{d}{dt}\int_{\mathbb{R}^2}\rho K_0(r)
dx\Big|_{t = 0} = \int_{\mathbb{R}^2}\rho_0 v_0K_0^\prime(r) dx>0.
\end{equation}
The combination of \eqref{4-4} and \eqref{4-5} gives that
\begin{eqnarray}\nonumber
&&G^\prime(t) = G^\prime(0) +  \int_0^t \frac{d^2}{ds^2}G(s) ds
\geq G^\prime(0) > 0.
\end{eqnarray}
Consequently, one can multiply the both sides of \eqref{4-4} by
$\frac{d}{dt}G(t)$ to get
\begin{eqnarray}\nonumber
&&\big(G^\prime(t)\big)^2\geq \frac{A}{(\gamma + 1)
  \big(\int_{\mathbb{R}^2}K_0(r)dx\big)^{\gamma-1}}\big(G(t)\big)^{\gamma + 1}\\\nonumber
&&\quad + \Big(\frac{d}{dt}\int_{\mathbb{R}^2}\rho K_0(r)
  dx\Big)^2\Big|_{t = 0} - \frac{A\big(
  \int_{\mathbb{R}^2}\rho_0K_0(r)
  dx\big)^{\gamma + 1}}{(\gamma + 1)
  \big(\int_{\mathbb{R}^2}K_0(r)dx\big)^{\gamma-1}}.
\end{eqnarray}
Noting \eqref{cond-6} and \eqref{4-5}, and by denoting
\begin{eqnarray}\nonumber
C_1 = \sqrt{\frac{A}{(\gamma + 1)
}}\big(\int_{\mathbb{R}^2}K_0(r)dx\big)^{- \frac{\gamma-1}{2}},
\end{eqnarray}
one has
\begin{eqnarray}\nonumber
G^\prime(t) \geq C_1\big(G(t)\big)^{\frac{\gamma + 1}{2}},
\end{eqnarray}
which implies that
\begin{equation}\nonumber
\int_{\mathbb{R}^2}\rho K_0(r)dx \geq \big(G(0)^{- \frac{\gamma -
1}{2}} - \frac{\gamma - 1}{2}C_1t\big)^{- \frac{2}{\gamma - 1}}.
\end{equation}
Noting \eqref{cond-5}, one has
\begin{equation}\nonumber
G(0) = \int_{\mathbb{R}^2} \rho_0K_0(r)dx > 0,
\end{equation}
By using the mass conservation, one has
\begin{equation}\nonumber
G(t) = \int_{\mathbb{R}^2}\frac{\rho}{re^{r}} dx \leq
\int_{B_{r_0}}\rho K_0(r) dx + \frac{\max_{r \geq
r_0}K_0(r)}{r_0}\int_{\mathbb{R}^2}\rho_0 dx,
\end{equation}
for any given $r_0 > 0$, here $B_{r_0}$ is the two-dimensional
ball centered at the origin with radius $r_0$. Consequently, one
concludes that $\rho(t, r)$ will blow up for $r \leq r_0$ as $t
\rightarrow \frac{2G(0)^{-\frac{\gamma - 1}{2}}}{(\gamma -
1)C_1}$. By the assumption $c \in L^\infty([0, \infty),
H^3(\mathbb{R}^2)$ and Sobolev imbedding theorem, one has $c \in
L^\infty([0, \infty)\times \mathbb{R}^2)$ and hence $\rho \in
L^\infty([0, \infty)\times \mathbb{R}^2)$. We arrive at a
contradiction. So the $H^3$ norm of $(c, u)$ can not be bounded
before the time $\frac{2G(0)^{-\frac{\gamma - 1}{2}}}{(\gamma -
1)C_1}$. The proof of Theorem \ref{Twod} is completed.

\end{proof}

\section*{Acknowledgement}
This work was done when Yi Du and Zhen Lei were visiting the
Department of Mathematics of Penn State University during 2012.
They would like to thank professor Qiang Du, professor Chun Liu
and the institute for their hospitality. Yi Du was supported by
NSFC (grant No. 11001088) and Pearl River New Star (grant No.
2012001). Zhen Lei was supported by NSFC (grant No.11171072), the
Foundation for Innovative Research Groups of NSFC (grant
No.11121101), FANEDD, Innovation Program of Shanghai Municipal
Education Commission (grant No.12ZZ012) and SGST 09DZ2272900.


\end{document}